\documentclass[12pt]{article}
\author{Edward Hanson}
\title{\bf A characterization of Leonard pairs\\ using the parameters $\{a_{i}\}_{i=0}^{d}$}
\date{}
\usepackage{amssymb}
\usepackage{amsmath}
\usepackage{cite}

\newtheorem{definition}{Definition}[section]
\newtheorem{theorem}[definition]{Theorem}
\newtheorem{proposition}[definition]{Proposition}
\newtheorem{lemma}[definition]{Lemma}
\newtheorem{corollary}[definition]{Corollary}

\newtheorem{example}[definition]{Example}

\newtheorem{note}[definition]{Note}
\newtheorem{assumption}[definition]{Assumption}

\usepackage[margin=1.0in]{geometry}

\def\fld{\mathbb K}

\begin{document}
\maketitle

\begin{abstract}
\noindent
Let $V$ denote a vector space with finite positive dimension. We consider an ordered pair of linear transformations 
$A: V\rightarrow V$ and $A^{*}: V\rightarrow V$ that satisfy (i) and (ii) below.
\begin{enumerate}
\item There exists a basis for $V$ with respect to which the matrix representing $A$ is irreducible tridiagonal and the matrix representing $A^{*}$ is diagonal.
\item There exists a basis for $V$ with respect to which the matrix representing $A^{*}$ is irreducible tridiagonal and the matrix representing $A$ is diagonal.
\end{enumerate}
We call such a pair a {\it Leonard pair} on $V$. Arlene Pascasio recently obtained a characterization of the $Q$-polynomial distance-regular graphs using the intersection numbers $a_{i}$. In this paper, we extend her results to a linear algebraic level and obtain a characterization of Leonard pairs. Pascasio's argument appears to rely on the underlying combinatorial assumptions, so we take a different approach that is algebraic in nature.

\bigskip

\end{abstract}

\section{Introduction} \label{sec:intro}

We begin by recalling the notion of a Leonard pair \cite{T:subconst1, T:Leonard}. We will use the following terms. Let $X$ denote a square matrix. Then $X$ is called {\it tridiagonal} whenever each nonzero entry lies on either the diagonal, the subdiagonal, or the superdiagonal. Assume $X$ is tridiagonal. Then $X$ is called {\it irreducible} whenever each entry on
the subdiagonal is nonzero and each entry on the superdiagonal is nonzero.

\medskip

\noindent We now define a Leonard pair. For the rest of this paper, $\fld$ will denote a field.

\begin{definition} \label{def:lp} \rm \cite[Definition 1.1]{T:Leonard}
Let $V$ denote a vector space over $\fld$ with finite positive dimension. By a {\it Leonard pair} on $V$, we mean an ordered pair of linear transformations $A: V\rightarrow V$ and $A^{*}: V\rightarrow V$ that satisfy (i) and (ii) below.
\begin{enumerate}
\item There exists a basis for $V$ with respect to which the matrix representing $A$ is irreducible tridiagonal and the matrix representing $A^{*}$ is diagonal.
\item There exists a basis for $V$ with respect to which the matrix representing $A^{*}$ is irreducible tridiagonal and the matrix representing $A$ is diagonal.
\end{enumerate}
\end{definition}

\begin{note}
\rm
It is a common notational convention to use $A^{*}$ to represent the conjugate-transpose of $A$. We are not using this convention. In a Leonard pair $A,A^{*}$, the linear transformations $A$ and $A^{*}$ are arbitrary subject to (i), (ii) above.
\end{note}

\noindent Leonard pairs appear in many contexts, such as orthogonal polynomials \cite{T:Leonard, T:qRacah}, Lie algebras \cite{Hartwig, T:intro}, quantum algebras \cite{IT}, and distance-regular graphs \cite{T:subconst1}. For a general survey, see \cite{T:madrid}. As researchers investigate Leonard pairs in these contexts, they devise characterization theorems for these objects that arise naturally in that context. There are characterizations of Leonard pairs in terms of orthogonal polynomials \cite[Theorem 19.1]{T:qRacah} \cite[Theorem 4.1]{T:PA}, the Lie algebra $\mathfrak{sl}_{2}$ \cite[Theorem 8.5]{Hartwig}, parameter arrays \cite[Theorem 1.9]{T:Leonard}, upper/lower bidiagonal matrices \cite[Theorem 3.2]{T:PA} \cite[Theorem 17.1]{T:TD-D}, tridiagonal/diagonal matrices \cite[Theorem 25.1]{T:TD-D}, and the notion of a tail \cite[Theorem 5.1]{Hanson}. 

\medskip

\noindent In the present paper, we obtain a characterization of Leonard pairs that is motivated by algebraic graph theory. Our result generalizes a result of Pascasio about $Q$-polynomial distance-regular graphs \cite[Theorem 1.2]{Pascasio}. In order to motivate our theorem, we first summarize Pascasio's result. Let $\Gamma$ denote a distance-regular graph with diameter $d\geq 3$ (see \cite{Pascasio} for definitions). Let $\theta$ denote a nontrivial eigenvalue of $\Gamma$ and let $\{\theta^{*}_{i}\}_{i=0}^{d}$ denote the corresponding dual eigenvalue sequence. Then by \cite[Theorem 1.2]{Pascasio}, $\Gamma$ is $Q$-polynomial with respect to $\theta$ if and only if the following (i)--(iii) hold.
\begin{enumerate}
\item[\rm (i)] There exist $\beta, \gamma^{*}\in {\mathbb C}$ such that 
\begin{equation} \label{eq:TTR_intro}
\gamma^{*}=\theta^{*}_{i-1}-\beta\theta^{*}_{i}+\theta^{*}_{i+1} \qquad \qquad (1\leq i\leq d-1).
\end{equation}
\item[\rm (ii)] There exist $\gamma, \omega, \eta^{*}\in {\mathbb C}$ such that the intersection numbers $a_{i}$ satisfy
\begin{equation}
a_{i}(\theta^{*}_{i}-\theta^{*}_{i-1})(\theta^{*}_{i}-\theta^{*}_{i+1})=\gamma \theta^{*2}_{i}+\omega \theta^{*}_{i}+\eta^{*} \qquad \qquad (0\leq i\leq d), \notag
\end{equation}
where $\theta^{*}_{-1}$ (resp. $\theta^{*}_{d+1}$) is the scalar which satisfies (\ref{eq:TTR_intro}) for $i=0$ (resp. $i=d$).
\item[\rm (iii)] $\theta^{*}_{i}\neq \theta^{*}_{0}$ for $1\leq i\leq d$.
\end{enumerate}
Very roughly speaking, a Leonard pair is a linear algebraic abstraction of a $Q$-polynomial distance-regular graph \cite[p.~260]{BIbook} \cite[Definition 2.3]{T:subconst1}. In the present paper, we obtain a characterization of Leonard pairs that is analogous to \cite[Theorem 1.2]{Pascasio}, but makes no reference to distance-regular graphs and is purely algebraic in nature. Our main result is Theorem \ref{thm:main} below. Note that Theorem \ref{thm:main} refers to the notion of a leaf. In the sections that follow that theorem, we obtain some results that are intended to illuminate the algebraic nature of these leaves.

\section{Leonard systems} \label{sec:LS}

\medskip
\noindent When working with a Leonard pair, it is often convenient to consider a closely related object called a Leonard system. To prepare for our definition of a Leonard system, we recall a few concepts from linear algebra. From now on, we fix a nonnegative integer $d$. Let $\mbox{Mat}_{d+1}(\fld)$ denote the $\fld$-algebra consisting of all $d+1$ by $d+1$ matrices with entries in $\fld$. We index the rows and columns by $0,1,\ldots ,d$. Let $\fld^{d+1}$ denote the $\fld$-vector space consisting of all $d+1$ by $1$ matrices with entries in $\fld$. We index the rows by $0,1,\ldots ,d$. Recall that $\mbox{Mat}_{d+1}(\fld)$ acts on $\fld^{d+1}$ by left multiplication. Let $V$ denote a vector space over $\fld$ with dimension $d+1$. Let $\mbox{End}(V)$ denote the $\fld$-algebra consisting of all linear transformations from $V$ to $V$. For convenience, we abbreviate $\mathcal{A}=\mbox{End}(V)$. Observe that $\mathcal{A}$ is $\fld$-algebra isomorphic to $\mbox{Mat}_{d+1}(\fld)$ and that $V$ is irreducible as an $\mathcal{A}$-module. The identity of $\mathcal{A}$ will be denoted by $I$. Let $\{v_{i}\}_{i=0}^{d}$ denote a basis for $V$. For $X\in \mathcal{A}$ and $Y\in \mbox{Mat}_{d+1}(\fld)$, we say that $Y$ {\it represents} $X$ {\it with respect to} $\{v_{i}\}_{i=0}^{d}$ whenever $Xv_{j}=\sum_{i=0}^{d}Y_{ij}v_{i}$ for $0\leq j\leq d$. Let $A$ denote an element of $\mathcal{A}$. A subspace $W\subseteq V$ will be called an {\it eigenspace} of $A$ whenever $W\neq 0$ and there exists $\theta \in \fld$ such that $W=\{v\in V|Av=\theta v\}$; in this case, $\theta$ is the {\it eigenvalue} of $A$ associated with $W$. We say that $A$ is {\it diagonalizable} whenever $V$ is spanned by the eigenspaces of $A$. We say that $A$ is {\it multiplicity-free} whenever it has $d+1$ mutually distinct eigenvalues in $\fld$. Note that if $A$ is multiplicity-free, then $A$ is diagonalizable.

\begin{definition} \label{def:MOidem}
\rm
By a {\it system of mutually orthogonal idempotents} in $\mathcal{A}$, we mean a sequence $\{E_{i}\}_{i=0}^{d}$ of elements in $\mathcal{A}$ such that
\begin{equation}
E_{i}E_{j}=\delta_{i,j}E_{i}\qquad \qquad (0 \leq i,j \leq d), \notag
\end{equation}
\begin{equation}
{\rm rank}(E_{i})=1\qquad \qquad (0 \leq i \leq d). \notag
\end{equation}
\end{definition}

\begin{definition} \label{def:decomp}
\rm
By a {\it decomposition of $V$}, we mean a sequence $\{U_{i}\}_{i=0}^{d}$ consisting of one-dimensional subspaces of $V$ such that
\begin{equation}
V=\sum_{i=0}^{d}U_{i}\qquad \qquad \text{(direct sum)}. \notag
\end{equation}
\end{definition}

\noindent The following lemmas are routinely verified.

\begin{lemma} \label{lem:MOidem2}
Let $\{U_{i}\}_{i=0}^{d}$ denote a decomposition of $V$. For $0\leq i\leq d$, define $E_{i}\in \mathcal{A}$ such that $(E_{i}-I)U_{i}=0$ and $E_{i}U_{j}=0$ if $j\ne i$ $(0\leq j\leq d)$. Then $\{E_{i}\}_{i=0}^{d}$ is a system of mutually orthogonal idempotents. Conversely, given a system of mutually orthogonal idempotents $\{E_{i}\}_{i=0}^{d}$ in $\mathcal{A}$, define $U_{i}=E_{i}V$ for $0\leq i\leq d$. Then $\{U_{i}\}_{i=0}^{d}$ is a decomposition of $V$.
\end{lemma}

\begin{lemma} \label{lem:EsumI}
Let $\{E_{i}\}_{i=0}^{d}$ denote a system of mutually orthogonal idempotents in $\mathcal{A}$. Then $I=\sum_{i=0}^{d}E_{i}$.
\end{lemma}

\noindent Let $A$ denote a multiplicity-free element of $\mathcal{A}$ and let $\{\theta_{i}\}^{d}_{i=0}$ denote an ordering of the eigenvalues of $A$. For $0\leq i\leq d$, let $U_{i}$ denote the eigenspace of $A$ for $\theta_{i}$. Then $\{U_{i}\}_{i=0}^{d}$ is a decomposition of $V$; let $\{E_{i}\}_{i=0}^{d}$ denote the corresponding system of idempotents from Lemma \ref{lem:MOidem2}. One checks that $A=\sum_{i=0}^{d}\theta_{i}E_{i}$ and $AE_{i}=E_{i}A=\theta_{i}E_{i}$ for $0\leq i\leq d$. Moreover,
\begin{equation}
E_{i}=\prod_{\genfrac{}{}{0pt}{}{0 \leq  j \leq d}{j\not=i}}\frac{A-\theta_{j}I}{\theta_{i}-\theta_{j}}\qquad \qquad (0\leq i\leq d). \label{eq:EpolyA}
\end{equation}
We refer to $E_{i}$ as the {\it primitive idempotent} of $A$ corresponding to $U_{i}$ (or $\theta_{i}$).

\medskip

\noindent We now define a Leonard system.

\begin{definition} \label{def:ls} \rm \cite[Definition 1.4]{T:Leonard}
By a {\it Leonard system} on $V$, we mean a sequence 
\begin{equation}
(A; \{E_{i}\}_{i=0}^{d}; A^{*}; \{E^{*}_{i}\}_{i=0}^{d}) \notag
\end{equation}
which satisfies (i)--(v) below.
\begin{enumerate}
\item Each of $A,A^{*}$ is a multiplicity-free element of $\mathcal{A}$.
\item $\{E_{i}\}_{i=0}^{d}$ is an ordering of the primitive idempotents of $A$.
\item $\{E^{*}_{i}\}_{i=0}^{d}$ is an ordering of the primitive idempotents of $A^{*}$.
\item ${\displaystyle{E^{*}_{i}AE^{*}_{j} =
\begin{cases}
0, & \text{if $\;|i-j|>1$;} \\
\neq 0, & \text{if $\;|i-j|=1$}
\end{cases}
}}
\qquad \qquad (0 \leq i,j\leq d).$
\item ${\displaystyle{E_{i}A^{*}E_{j} =
\begin{cases}
0, & \text{if $\;|i-j|>1$;} \\
\neq 0, & \text{if $\;|i-j|=1$}
\end{cases}
}}
\qquad \qquad (0 \leq i,j\leq d).$
\end{enumerate}
\end{definition}

\noindent Leonard systems and Leonard pairs are related as follows. Let $(A; \{E_{i}\}_{i=0}^{d}; A^{*}; \{E^{*}_{i}\}_{i=0}^{d})$ denote a Leonard system on $V$. For $0\leq i\leq d$, let $v_{i}$ denote a nonzero vector in $E_{i}V$. Then the sequence $\{v_{i}\}_{i=0}^{d}$ is a basis for $V$ which satisfies Definition \ref{def:lp}(ii). For $0\leq i\leq d$, let $v^{*}_{i}$ denote a nonzero vector in $E^{*}_{i}V$. Then the sequence $\{v^{*}_{i}\}_{i=0}^{d}$ is a basis for $V$ which satisfies Definition \ref{def:lp}(i). By these comments, the pair $A,A^{*}$ is a Leonard pair on $V$. Conversely, let $A,A^{*}$ denote a Leonard pair on $V$. By \cite[Lemma 1.3]{T:Leonard}, each of $A,A^{*}$ is multiplicity-free. Let $\{v_{i}\}_{i=0}^{d}$ denote a basis for $V$ which satisfies Definition \ref{def:lp}(ii). For $0\leq i\leq d$, the vector $v_{i}$ is an eigenvector for $A$; let $E_{i}$ denote the corresponding primitive idempotent. Let $\{v^{*}_{i}\}_{i=0}^{d}$ denote a basis for $V$ which satisfies Definition \ref{def:lp}(i). For $0\leq i\leq d$, the vector $v^{*}_{i}$ is an eigenvector for $A^{*}$; let $E^{*}_{i}$ denote the corresponding primitive idempotent. Then $(A; \{E_{i}\}_{i=0}^{d}; A^{*}; \{E^{*}_{i}\}_{i=0}^{d})$ is a Leonard system on $V$.

\medskip

\noindent We make some observations. Let $(A; \{E_{i}\}_{i=0}^{d}; A^{*}; \{E^{*}_{i}\}_{i=0}^{d})$ denote a Leonard system on $V$. For $0\leq i\leq d$, let $\theta_{i}$ (resp. $\theta^{*}_{i}$) denote the eigenvalue of $A$ (resp. $A^{*}$) associated with $E_{i}V$ (resp. $E^{*}_{i}V$). By construction, $\{\theta_{i}\}_{i=0}^{d}$ (resp. $\{\theta^{*}_{i}\}_{i=0}^{d}$) are mutually distinct and contained in $\fld$. By \cite[Theorem 12.7]{T:Leonard}, the expressions
\begin{equation} \label{eq:thetarecur}
\frac{\theta_{i-2}-\theta_{i+1}}{\theta_{i-1}-\theta_{i}}, \qquad \frac{\theta^{*}_{i-2}-\theta^{*}_{i+1}}{\theta^{*}_{i-1}-\theta^{*}_{i}}
\end{equation}
are equal and independent of $i$ for $2\leq i\leq d-1$. Define $\beta\in\fld$ such that $\beta+1$ equals the common value of (\ref{eq:thetarecur}). If $d\leq 2$, let $\beta$ be arbitrary. It will be useful to describe the above features as follows. By \cite[Lemmas 8.3 and 8.4]{T:Leonard}, there exists $\gamma\in\fld$ such that
\begin{equation} \label{eq:gamma}
\theta_{i-1}-\beta \theta_{i}+\theta_{i+1}=\gamma \qquad \qquad (1\leq i\leq d-1)
\end{equation}
and there exists $\gamma^{*}\in\fld$ such that
\begin{equation} \label{eq:introTTR}
\theta^{*}_{i-1}-\beta \theta^{*}_{i}+\theta^{*}_{i+1}=\gamma^{*} \qquad \qquad (1\leq i\leq d-1).
\end{equation}

\section{The antiautomorphism $\dagger$}

In this section, we discuss an antiautomorphism related to Leonard systems.

\begin{lemma} \label{lem:TDpower} {\rm \cite[Lemma 5.3]{T:madrid}}
Let $A$ denote an irreducible tridiagonal matrix in $\mbox{\rm{Mat}}_{d+1}(\fld)$. Then the following {\rm (i)--(iii)} hold for $0\leq i,j\leq d$.
\begin{enumerate}
\item[\rm (i)] The entry $(A^{r})_{i,j}=0$ if $r<|i-j|$\qquad \qquad $(0 \leq r \leq d)$.
\item[\rm (ii)] Suppose $i\leq j$. Then the entry $(A^{j-i})_{i,j}=\prod_{h=i}^{j-1}A_{h,h+1}$. Moreover, $(A^{j-i})_{i,j} \ne 0$.
\item[\rm (iii)] Suppose $i\geq j$. Then the entry $(A^{i-j})_{i,j}=\prod_{h=j}^{i-1}A_{h+1,h}$. Moreover, $(A^{i-j})_{i,j} \ne 0$.
\end{enumerate}
\end{lemma}
\noindent {\it Proof:} This follows from the definition of matrix multiplication and the meaning of irreducible tridiagonal. \hfill $\Box$ \\

\begin{assumption} \label{assum:E*AE*}
\rm
Let $\{E^{*}_{i}\}^{d}_{i=0}$ denote a system of mutually orthogonal idempotents in $\mathcal{A}$. Let $A$ denote an element of $\mathcal{A}$ such that
\begin{equation} \label{eq:E*AE*}
{\displaystyle{E^{*}_{i}AE^{*}_{j} =
\begin{cases}
0, & \text{if $\;|i-j|>1$;} \\
\neq 0, & \text{if $\;|i-j|=1$}
\end{cases}
}}
\qquad \qquad (0 \leq i,j\leq d).
\end{equation}
\end{assumption}

\noindent We make some comments on Assumption \ref{assum:E*AE*}. Let $\{E^{*}_{i}\}_{i=0}^{d}$ and $A$ denote elements of $\mathcal{A}$ that satisfy Assumption \ref{assum:E*AE*}. For $0\leq j\leq d$, let $v^{*}_{j}$ denote a nonzero vector in $E^{*}_{j}V$ and note that $\{v^{*}_{j}\}_{j=0}^{d}$ is a basis for $V$. For $0\leq i\leq d$, the matrix representing $E^{*}_{i}$ with respect to $\{v^{*}_{j}\}_{j=0}^{d}$ is diagonal, with $(i,i)$-entry $1$ and all other entries $0$. By this and (\ref{eq:E*AE*}), the matrix representing $A$ with respect to $\{v^{*}_{j}\}_{j=0}^{d}$ is irreducible tridiagonal.

\medskip

\noindent Conversely, let $\{v^{*}_{j}\}_{j=0}^{d}$ denote a basis for $V$. For $0\leq i\leq d$, define $E^{*}_{i}\in \mathcal{A}$ such that $E^{*}_{i}v^{*}_{i}=v^{*}_{i}$ and $E^{*}_{i}v^{*}_{j}=0$ if $j\neq i$ ($0\leq j\leq d$). The matrix representing $E^{*}_{i}$ with respect to $\{v^{*}_{j}\}_{j=0}^{d}$ is diagonal, with $(i,i)$-entry $1$ and all other entries $0$. The sequence $\{E^{*}_{i}\}_{i=0}^{d}$ satisfies Definition \ref{def:MOidem}, so it is a system of mutually orthogonal idempotents. Let $A\in \mathcal{A}$ denote the linear transformation represented by an irreducible tridiagonal matrix with respect to the basis $\{v^{*}_{j}\}_{j=0}^{d}$. Then $\{E^{*}_{i}\}_{i=0}^{d}$ and $A$ satisfy (\ref{eq:E*AE*}). By these comments, $\{E^{*}_{i}\}_{i=0}^{d}$ and $A$ satisfy Assumption \ref{assum:E*AE*}.

\begin{lemma} \label{lem:AE*0gen} {\rm \cite[Corollary 3.4]{Hanson}}
With reference to Assumption \ref{assum:E*AE*}, the elements $A$ and $E^{*}_{0}$ together generate $\mathcal{A}$.
\end{lemma}

\noindent We recall the notion of an {\it antiautomorphism} of $\mathcal{A}$. Let $\xi : \mathcal{A} \rightarrow \mathcal{A}$ denote any map. We call $\xi$ an {\it antiautomorphism} of $\mathcal{A}$ whenever $\xi$ is an isomorphism of $\fld$-vector spaces and $(XY)^{\xi}=Y^{\xi}X^{\xi}$ for all $X, Y\in \mathcal{A}$.

\begin{lemma} \label{lem:dagexist} {\rm \cite[Lemma 3.5]{Hanson}}
With reference to Assumption \ref{assum:E*AE*}, there exists a unique antiautomorphism $\dagger$ of $\mathcal{A}$ such that $A^{\dagger}=A$ and $E^{*\dagger}_{0}=E^{*}_{0}$. Moreover, $E^{*\dagger}_{i}=E^{*}_{i}$ for $1 \leq i \leq d$ and $X^{\dagger \dagger}=X$ for all $X \in \mathcal{A}$. 
\end{lemma}

\begin{definition} \label{def:a}
\rm
With reference to Assumption \ref{assum:E*AE*}, define
\begin{equation}
a_{i}=\mbox{tr}(E^{*}_{i}A) \qquad \qquad (0 \leq i \leq d), \notag
\end{equation}
where $\mbox{tr}$ denotes trace.
\end{definition}

\begin{proposition} \label{prop:E*AE*}
With reference to Assumption \ref{assum:E*AE*}, $E^{*}_{i}AE^{*}_{i}=a_{i}E^{*}_{i}$ for $0\leq i\leq d$.
\end{proposition}
\noindent {\it Proof:} For $0\leq i\leq d$, let $v^{*}_{i}$ denote a nonzero vector in $E^{*}_{i}V$ and note that $\{v^{*}_{i}\}_{i=0}^{d}$ is a basis for $V$. By the discussion following Assumption \ref{assum:E*AE*}, for $0\leq i\leq d$, the matrix representing $E^{*}_{i}$ with respect to this basis is diagonal with $(i,i)$-entry $1$ and all other entries $0$. Using matrix multiplication, we observe that $E^{*}_{i}AE^{*}_{i} = \alpha E^{*}_{i}$ for some $\alpha \in \fld$. Taking the trace of both sides establishes the result. \hfill $\Box$ \\

\noindent We have been discussing the situation of Assumption \ref{assum:E*AE*}. We now modify this situation as follows.

\begin{assumption} \label{assum:A*}
\rm
Let $A$ and $\{E^{*}_{i}\}_{i=0}^{d}$ be as in Assumption \ref{assum:E*AE*}. Furthermore, assume that $A$ is multiplicity-free, with primitive idempotents $\{E_{i}\}_{i=0}^{d}$ and eigenvalues $\{\theta_{i}\}_{i=0}^{d}$. Additionally, let $\{\theta^{*}_{i}\}^{d}_{i=0}$ denote scalars in $\fld$ and let $A^{*}=\sum_{i=0}^{d}\theta^{*}_{i}E^{*}_{i}$. To avoid trivialities, assume that $d\geq 1$.
\end{assumption}

\begin{lemma} \label{lem:eistab} {\rm \cite[Lemma 3.7]{Hanson}}
With reference to Assumption \ref{assum:A*}, the antiautomorphism $\dagger$ from Lemma \ref{lem:dagexist} satisfies $A^{*\dagger}=A^{*}$ and $E_{i}^{\dagger}=E_{i}$ for $0 \leq i \leq d$.
\end{lemma}
\noindent {\it Proof:} By (\ref{eq:EpolyA}), $E_{i}$ is a polynomial in $A$ for $0 \leq i \leq d$. The result follows in view of Lemma \ref{lem:dagexist}. \hfill $\Box$ \\

\begin{lemma} \label{lem:undirected} {\rm \cite[Lemma 3.8]{Hanson}}
With reference to Assumption \ref{assum:A*} and for $0\leq i,j\leq d$, $E_{i}A^{*}E_{j}=0$ if and only if $E_{j}A^{*}E_{i}=0$.
\end{lemma}
\noindent {\it Proof:} Let $\dagger$ denote the antiautomorphism from Lemma \ref{lem:dagexist}. Then $E_{i}A^{*}E_{j}=0$ if and only if $(E_{i}A^{*}E_{j})^{\dagger}=0$. Also, using Lemma \ref{lem:eistab}, $(E_{i}A^{*}E_{j})^{\dagger}=E_{j}^{\dagger}A^{* \dagger}E_{i}^{\dagger}=E_{j}A^{*}E_{i}$. The result follows. \hfill $\Box$ \\

\noindent We mention a result for later use.

\begin{lemma} \label{lem:TV} {\rm \cite[Theorem 5.3 (iv)]{TV}}
With reference to Assumption \ref{assum:A*}, further assume that the sequence $(A; \{E_{i}\}^{d}_{i=0}; A^{*}; \{E^{*}_{i}\}^{d}_{i=0})$ is a Leonard system. Then there exist $\omega, \eta^{*} \in \fld$ such that
\begin{equation}
a_{i}(\theta^{*}_{i}-\theta^{*}_{i-1})(\theta^{*}_{i}-\theta^{*}_{i+1})=\gamma \theta^{*2}_{i}+\omega \theta^{*}_{i}+\eta^{*} \qquad \qquad (0\leq i\leq d), \notag
\end{equation}
where $\gamma$ is from {\rm (\ref{eq:gamma})}, and $\theta^{*}_{-1}$ (resp. $\theta^{*}_{d+1}$) is the scalar which satisfies {\rm (\ref{eq:introTTR})} for $i=0$ (resp. $i=d$).
\end{lemma}

\section{The graph $\Delta$} \label{sec:delta}

In the following discussion, a graph is understood to be finite and undirected, without loops or multiple edges.

\begin{definition} \label{def:delta}
\rm
With reference to Assumption \ref{assum:A*}, let $\Delta$ denote the graph with vertex set $\{0, 1,\ldots ,d\}$ such that two vertices $i, j$ are adjacent if and only if $i\neq j$ and $E_{i}A^{*}E_{j}\neq 0$. The graph $\Delta$ is well-defined in view of Lemma \ref{lem:undirected}.
\end{definition}

\begin{lemma} \label{lem:ijadjacent}
With reference to Assumption \ref{assum:A*}, let $i$ and $j$ denote distinct vertices in $\Delta$. Then $i, j$ are adjacent if and only if $E_{i}A^{*}E_{j}V=E_{i}V$.
\end{lemma}
\noindent {\it Proof:} Suppose $i$ and $j$ are adjacent. Then $E_{i}A^{*}E_{j}\neq 0$ by Definition \ref{def:delta}. Now $E_{i}A^{*}E_{j}V$ is a nonzero subspace of the one-dimensional space $E_{i}V$, so $E_{i}A^{*}E_{j}V=E_{i}V$.

\medskip

\noindent Conversely, assume that $E_{i}A^{*}E_{j}V=E_{i}V$. Then $E_{i}A^{*}E_{j}\neq 0$, so $i$ and $j$ are adjacent by Definition \ref{def:delta}. \hfill $\Box$ \\

\begin{lemma} \label{lem:lspath} {\rm \cite[Lemma 4.2]{Hanson}}
With reference to Assumption \ref{assum:A*}, the following {\rm (i)}, {\rm (ii)} are equivalent.
\begin{enumerate}
\item[\rm (i)] The sequence $(A; \{E_{i}\}^{d}_{i=0}; A^{*}; \{E^{*}_{i}\}^{d}_{i=0})$ is a Leonard system.
\item[\rm (ii)] The graph $\Delta$ is a path such that vertices $i-1, i$ are adjacent for $1\leq i\leq d$.
\end{enumerate}
\end{lemma}

\begin{definition} \label{def:Qpoly1}
\rm
With reference to Assumption \ref{assum:A*}, $A^{*}$ is said to be {\it $Q$-polynomial} whenever $\Delta$ is a path.
\end{definition}

\begin{definition} \label{def:stail}
\rm
With reference to Assumption \ref{assum:A*}, let $E=E_{i}$ denote a primitive idempotent for $A$. This idempotent will be called a {\it leaf} whenever $i$ is adjacent to at most one vertex in $\Delta$.
\end{definition}

\begin{example} \label{ex:Qpolystail}
\rm
With reference to Assumption \ref{assum:A*}, assume further that $A^{*}$ is $Q$-polynomial. By Definition \ref{def:Qpoly1}, $\Delta$ is a path. Fix an endpoint of the path $\Delta$ and relabel the primitive idempotents of $A$ such that this endpoint is vertex $0$ and vertices $i-1, i$ are adjacent for $1\leq i\leq d$. By Lemma \ref{lem:lspath}, the sequence $(A; \{E_{i}\}^{d}_{i=0}; A^{*}; \{E^{*}_{i}\}^{d}_{i=0})$ is a Leonard system. Also note that both $E_{0}$ and $E_{d}$ are leaves.
\end{example}

\noindent For the rest of this section, we discuss the connectivity of $\Delta$.

\begin{lemma} \label{lem:invconnect1} {\rm \cite[Lemma 4.7]{Hanson}}
With reference to Assumption \ref{assum:A*}, fix a subspace $U\subseteq V$. Then $AU \subseteq U$ if and only if there exists a subset $S\subseteq \{0, 1,\ldots ,d\}$ such that $U=\sum_{h\in S}E_{h}V$. In this case, $S$ is uniquely determined by $U$.
\end{lemma}

\noindent We will use the following notation. For a subset $S\subseteq \{0, 1,\ldots ,d\}$, let $\overline{S}$ denote the complement of $S$ in $\{0, 1,\ldots ,d\}$.

\begin{proposition} \label{prop:invconnect2} {\rm \cite[Proposition 4.8]{Hanson}}
With reference to Assumption \ref{assum:A*}, fix a subset $S\subseteq \{0, 1,\ldots ,d\}$ and let $U=\sum_{h\in S}E_{h}V$. Then the following {\rm (i)}, {\rm (ii)} are equivalent.
\begin{enumerate}
\item[\rm (i)] $A^{*}U \subseteq U$.
\item[\rm (ii)] The vertices $i,j$ are not adjacent in the graph $\Delta$ for all $i\in S$ and $j\in \overline{S}$.
\end{enumerate}
\end{proposition}

\begin{proposition} \label{prop:deltaconnected}
With reference to Assumption \ref{assum:A*}, assume further that $\theta^{*}_{i}\neq \theta^{*}_{0}$ for $1\leq i \leq d$. Then $\Delta$ is connected.
\end{proposition}
\noindent {\it Proof:} Suppose $\Delta$ is not connected. Then there exists a non-empty proper subset $S$ of $\{0,1,\ldots ,d\}$ such that $i, j$ are not adjacent in $\Delta$ for all $i\in S$ and $j\in \overline{S}$. Let $U=\sum_{h\in S}E_{h}V$ and note that $U\neq 0$ and $U\neq V$. Observe that $AU\subseteq U$ by Lemma \ref{lem:invconnect1} and $A^{*}U\subseteq U$ by Proposition \ref{prop:invconnect2}. Using the equation $A^{*}=\sum_{i=0}^{d}\theta^{*}_{i}E^{*}_{i}$ and the fact that $\{E^{*}_{i}\}_{i=0}^{d}$ are mutually orthogonal idempotents,
\begin{equation}
E^{*}_{0}=\prod_{j=1}^{d}\frac{A^{*}-\theta^{*}_{j}I}{\theta^{*}_{0}-\theta^{*}_{j}}. \label{eq:E*0polyA*}
\end{equation}
The denominator in (\ref{eq:E*0polyA*}) is nonzero by assumption. By (\ref{eq:E*0polyA*}) and since $A^{*}U\subseteq U$, we find that $E^{*}_{0}U\subseteq U$. By Lemma \ref{lem:AE*0gen}, $A$ and $E^{*}_{0}$ generate $\mathcal{A}$. Therefore, $\mathcal{A}U\subseteq U$. Recall that $V$ is irreducible as an $\mathcal{A}$-module, so either $U=0$ or $U=V$. This is a contradiction, so $\Delta$ is connected. \hfill $\Box$ \\

\section{The main theorem}

The following is our main result.

\begin{theorem} \label{thm:main}
With reference to Assumption \ref{assum:A*}, $A^{*}$ is $Q$-polynomial if and only if the following {\rm (i)--(iv)} hold.
\begin{enumerate}
\item[\rm (i)] There exists a leaf in $\Delta$.
\item[\rm (ii)] There exist $\beta, \gamma^{*} \in \fld$ such that
\begin{equation} \label{eq:TTR}
\gamma^{*}=\theta^{*}_{i-1}-\beta \theta^{*}_{i}+\theta^{*}_{i+1} \qquad \qquad (1\leq i \leq d-1).
\end{equation}
\item[\rm (iii)] There exist $\gamma, \omega, \eta^{*} \in \fld$ such that
\begin{equation} \label{eq:ai_assum}
a_{i}(\theta^{*}_{i}-\theta^{*}_{i-1})(\theta^{*}_{i}-\theta^{*}_{i+1})=\gamma \theta^{*2}_{i}+\omega \theta^{*}_{i}+\eta^{*} \qquad \qquad (0\leq i\leq d),
\end{equation}
where $\theta^{*}_{-1}$ (resp. $\theta^{*}_{d+1}$) is the scalar which satisfies {\rm (\ref{eq:TTR})} for $i=0$ (resp. $i=d$).
\item[\rm (iv)] $\theta^{*}_{i}\neq \theta^{*}_{0}$ for $1\leq i \leq d$.
\end{enumerate}
\end{theorem}
\noindent {\it Proof:} First, assume that $A^{*}$ is $Q$-polynomial, so that $\Delta$ is a path. Label the vertex set of $\Delta$ such that vertices $i-1$ and $i$ are adjacent for $1\leq i\leq d$. Note that $E_{0}$ is a leaf, so condition (i) is satisfied. By Lemma \ref{lem:lspath}, the sequence $(A; \{E_{i}\}^{d}_{i=0}; A^{*}; \{E^{*}_{i}\}^{d}_{i=0})$ is a Leonard system. Now condition (ii) follows from (\ref{eq:introTTR}). We mentioned near the end of Section \ref{sec:LS} that $\{\theta^{*}_{i}\}_{i=0}^{d}$ are mutually distinct, from which (iv) follows. Condition (iii) follows from Lemma \ref{lem:TV}.

\medskip

\noindent Conversely, assume that conditions (i)--(iv) hold. We show that $A^{*}$ is $Q$-polynomial. To do this, we show that $\Delta$ is a path. First note that $\Delta$ is connected by Proposition \ref{prop:deltaconnected}.

\medskip

\noindent Define $\theta^{*}_{-1}$ and $\theta^{*}_{d+1}$ such that (\ref{eq:TTR}) holds for $i=0$ and $i=d$, so that 
\begin{equation} \label{eq:TTR2}
\gamma^{*}=\theta^{*}_{i-1}-\beta \theta^{*}_{i}+\theta^{*}_{i+1} \qquad \qquad (0\leq i \leq d).
\end{equation}
We claim that the expression
\begin{equation} \label{eq:delta*}
\theta ^{*2}_{i-1}-\beta \theta ^{*}_{i-1}\theta ^{*}_{i}+\theta ^{*2}_{i}-\gamma ^{*}(\theta ^{*}_{i-1}+\theta ^{*}_{i}) 
\end{equation}
is independent of $i$ for $0\leq i\leq d+1$. Denote this expression by $p_{i}$. Observe that, for $0\leq i\leq d$,
\begin{equation}
p_{i}-p_{i+1}=(\theta^{*}_{i-1}-\theta^{*}_{i+1})(\theta^{*}_{i-1}-\beta\theta^{*}_{i}+\theta^{*}_{i+1}-\gamma^{*}). \notag
\end{equation}
In this equation, the expression on the right-hand side equals $0$ by (\ref{eq:TTR2}). Consequently, $p_{i}$ is independent of $i$ for $0\leq i\leq d+1$. The claim is now proven. Let $\delta^{*}$ denote the common value of (\ref{eq:delta*}).

\medskip

\noindent We now show that
\begin{equation} \label{eq:theta*prod}
(\theta^{*}_{i}-\theta^{*}_{i-1})(\theta^{*}_{i}-\theta^{*}_{i+1})=(2-\beta)\theta^{*2}_{i}-2\gamma^{*}\theta^{*}_{i}-\delta^{*} \qquad \qquad (0\leq i \leq d).
\end{equation}
To verify (\ref{eq:theta*prod}), in the right-hand side, replace $\delta^{*}$ by (\ref{eq:delta*}) and eliminate both occurrences of $\gamma^{*}$ in the resulting expression using $\gamma^{*}=\theta^{*}_{i-1}-\beta\theta^{*}_{i}+\theta^{*}_{i+1}$. We have now verified (\ref{eq:theta*prod}).

\medskip

\noindent For notational convenience, we introduce a $2$-variable polynomial 
\begin{equation}
P(\lambda,\mu)=\lambda^{2}-\beta \lambda \mu+\mu^{2}-\gamma^{*}(\lambda +\mu)-\delta^{*}. \notag
\end{equation}
We now claim that
\begin{equation} \label{eq:AW2}
A^{*2}A-\beta A^{*}AA^{*}+AA^{*2}-\gamma^{*}(AA^{*}+A^{*}A)-\delta^{*}A=\gamma A^{*2}+\omega A^{*}+\eta^{*}I.
\end{equation}
In (\ref{eq:AW2}), let $C$ denote the left-hand side minus the right-hand side. We show $C=0$. Using $I=\sum_{i=0}^{d}E^{*}_{i}$, we obtain
\begin{align}
C &= (E^{*}_{0}+E^{*}_{1}+\cdots +E^{*}_{d})C(E^{*}_{0}+E^{*}_{1}+\cdots +E^{*}_{d}) \notag \\
  &= \sum_{i=0}^{d}\sum_{j=0}^{d} E^{*}_{i}CE^{*}_{j}. \notag
\end{align}
For $0\leq i,j\leq d$, we show $E^{*}_{i}CE^{*}_{j}=0$. Recall that $E^{*}_{i}A^{*}=\theta^{*}_{i}E^{*}_{i}$ and $A^{*}E^{*}_{j}=\theta^{*}_{j}E^{*}_{j}$. 
Thus,
\begin{equation} \label{eq:E*iCE*j}
E^{*}_{i}CE^{*}_{j}=E^{*}_{i}AE^{*}_{j} P(\theta^{*}_{i},\theta^{*}_{j})-\delta_{i,j}(\gamma \theta^{*2}_{i}+\omega \theta^{*}_{i}+\eta^{*})E^{*}_{i}.
\end{equation}
To further examine (\ref{eq:E*iCE*j}), we consider two cases. First assume $i\neq j$. In this case, (\ref{eq:E*iCE*j}) becomes
\begin{equation}
E^{*}_{i}CE^{*}_{j}=E^{*}_{i}AE^{*}_{j} P(\theta^{*}_{i},\theta^{*}_{j}). \notag
\end{equation}
If $|i-j|>1$, then $E^{*}_{i}AE^{*}_{j}=0$ by Assumption \ref{assum:A*}. If $|i-j|=1$, then $P(\theta^{*}_{i},\theta^{*}_{j})=0$. Therefore, $E^{*}_{i}CE^{*}_{j}=0$ under our present assumption that $i\neq j$. Next assume $i=j$. In this case, (\ref{eq:E*iCE*j}) becomes
\begin{equation} \label{eq:E*iCE*i}
E^{*}_{i}CE^{*}_{i}=E^{*}_{i}AE^{*}_{i}P(\theta^{*}_{i},\theta^{*}_{i})-(\gamma \theta^{*2}_{i}+\omega \theta^{*}_{i}+\eta^{*})E^{*}_{i}.
\end{equation}
By the definition of $P$ and (\ref{eq:theta*prod}), we find $P(\theta^{*}_{i},\theta^{*}_{i})=(\theta^{*}_{i}-\theta^{*}_{i-1})(\theta^{*}_{i}-\theta^{*}_{i+1})$. Using Proposition \ref{prop:E*AE*}, $E^{*}_{i}AE^{*}_{i}=a_{i}E^{*}_{i}$. Evaluating the right-hand side of (\ref{eq:E*iCE*i}) using these comments, we find that it equals $E^{*}_{i}$ times
\begin{equation} \label{eq:aithetas}
a_{i}(\theta^{*}_{i}-\theta^{*}_{i-1})(\theta^{*}_{i}-\theta^{*}_{i+1})-\gamma \theta^{*2}_{i}-\omega \theta^{*}_{i}-\eta^{*}.
\end{equation}
Note that (\ref{eq:aithetas}) is equal to 0 by (\ref{eq:ai_assum}). Therefore, $E^{*}_{i}CE^{*}_{i}=0$. We have now shown $E^{*}_{i}CE^{*}_{j}=0$ for $0\leq i,j\leq d$. Therefore, $C=0$. We have now verified (\ref{eq:AW2}).

\medskip

\noindent Suppose we are given vertices $i$ and $j$ in $\Delta$ with $\partial (i,j)=2$, where $\partial$ denotes path-length distance. Further, suppose there exists a unique vertex $r\in\Delta$ adjacent to both $i$ and $j$. We show
\begin{equation} \label{eq:irj}
\gamma=\theta_{i}-\beta\theta_{r}+\theta_{j}.
\end{equation}
To show (\ref{eq:irj}), we multiply (\ref{eq:AW2}) on the left by $E_{i}$ and on the right by $E_{j}$, and simplify. To illustrate, we simplify the second term. Using $A=\sum_{h=0}^{d}\theta_{h}E_{h}$ and Definition \ref{def:delta},
\begin{align}
E_{i}A^{*}AA^{*}E_{j} &= E_{i}A^{*}\left( \sum_{h=0}^{d}\theta_{h}E_{h} \right) A^{*}E_{j} \notag \\
                      &= \theta_{r}E_{i}A^{*}E_{r}A^{*}E_{j}. \notag
\end{align}
Simplifying the other terms in a similar fashion yields
\begin{align}
    E_{i}A^{*2}AE_{j} &= \theta_{j}E_{i}A^{*}E_{r}A^{*}E_{j}, \notag \\
    E_{i}AA^{*2}E_{j} &= \theta_{i}E_{i}A^{*}E_{r}A^{*}E_{j}, \notag \\
     E_{i}A^{*2}E_{j} &= E_{i}A^{*}E_{r}A^{*}E_{j}, \notag
\end{align}
\begin{align}
E_{i}AA^{*}E_{j} &= 0, & E_{i}A^{*}AE_{j}=0, \notag \\
     E_{i}AE_{j} &= 0, & E_{i}A^{*}E_{j}=0. \notag
\end{align}
By the above comments,
\begin{equation}
(\theta_{i}-\beta\theta_{r}+\theta_{j}-\gamma) E_{i}A^{*}E_{r}A^{*}E_{j}=0. \label{eq:AW2result}
\end{equation}
We now show $E_{i}A^{*}E_{r}A^{*}E_{j} \neq 0$. Since $r$ and $j$ are adjacent, $E_{r}A^{*}E_{j}V=E_{r}V$ by Lemma \ref{lem:ijadjacent}. Similarly, $E_{i}A^{*}E_{r}V=E_{i}V$, so $E_{i}A^{*}E_{r}A^{*}E_{j}V=E_{i}V$. Therefore, $E_{i}A^{*}E_{r}A^{*}E_{j}\neq 0$. This and (\ref{eq:AW2result}) imply (\ref{eq:irj}).

\medskip

\noindent By (i), there exists a leaf, which we call $E_{0}$. Because $E_{0}$ is a leaf and $\Delta$ is connected, vertex $0$ is adjacent to a single vertex. We can now easily show that $\Delta$ is a path. To this end, we show that every vertex in $\Delta$ is adjacent to at most two other vertices in $\Delta$. Suppose there exists a vertex $i$ in $\Delta$ that is adjacent to at least three other vertices. Choose the $i$ such that $\partial (0,i)$ is minimum. Without loss of generality, assume that the nonzero vertices of $\Delta$ are labeled such that $\partial (0,i)=i$ and $(0,1,\ldots,i)$ is a path. By construction, $i\geq 1$. By assumption, there exist distinct vertices $j$ and $j'$, each at least $i+1$, that are both adjacent to $i$. Note that $\partial (i-1,j)=2$ and that $i$ is the unique vertex in $\Delta$ adjacent to $i-1$ and $j$. Therefore, by (\ref{eq:irj}), 
\begin{equation}
\theta_{i-1}-\beta\theta_{i}+\theta_{j}=\gamma. \label{eq:thetaj}
\end{equation}
Replacing $j$ by $j'$ in the above argument, we obtain 
\begin{equation}
\theta_{i-1}-\beta\theta_{i}+\theta_{j'}=\gamma. \label{eq:thetaj'}
\end{equation}
Comparing (\ref{eq:thetaj}) to (\ref{eq:thetaj'}), we find $\theta _{j}=\theta _{j'}$. Recall that $\{\theta_{h}\}_{h=0}^{d}$ are mutually distinct, so $j=j'$. This is a contradiction and we have now shown that $\Delta$ is a path. Therefore, $A^{*}$ is $Q$-polynomial. \hfill $\Box$ \\

\section{Recognizing leaves in $\Delta$ (part 1)} \label{sec:a*}

We wish to gain a more thorough understanding of Theorem \ref{thm:main}(i). With reference to Assumption \ref{assum:A*}, suppose we are given two distinct vertices of $\Delta$, denoted $r$ and $s$. Our goal for the remainder of the paper is to develop necessary and sufficient conditions for $r$ to be adjacent to $s$ and no other vertices. We will examine this from several different perspectives.

\begin{definition} \label{def:a*}
\rm
With reference to Assumption \ref{assum:A*}, define
\begin{equation}
a^{*}_{i}=\mbox{tr}(E_{i}A^{*}) \qquad \qquad (0 \leq i \leq d). \notag
\end{equation}
\end{definition}

\begin{proposition} \label{prop:EA*E}
With reference to Assumption \ref{assum:A*}, $E_{i}A^{*}E_{i}=a^{*}_{i}E_{i}$ for $0\leq i\leq d$.
\end{proposition}
\noindent {\it Proof:} Similar to Proposition \ref{prop:E*AE*}. \hfill $\Box$ \\

\begin{lemma} \label{lem:A*-a*I}
With reference to Assumption \ref{assum:A*}, fix distinct integers $r$ and $s$ such that $0\leq r,s\leq d$. Then the following {\rm (i)}, {\rm (ii)} are equivalent.
\begin{enumerate}
\item[\rm (i)] In the diagram $\Delta$, vertex $r$ is adjacent to vertex $s$ and no other vertices.
\item[\rm (ii)] There exists $\kappa\in\fld$ such that $(A^{*}-\kappa I)E_{r}V=E_{s}V$.
\end{enumerate}
Suppose conditions {\rm (i)} and {\rm (ii)} hold. Then $\kappa=a^{*}_{r}$.
\end{lemma}
\noindent {\it Proof:} (i) $\Rightarrow$ (ii). Using $I=\sum_{i=0}^{d} E_{i}$ and Definition \ref{def:delta},
\begin{align}
(A^{*}-a^{*}_{r} I)E_{r} &= \sum_{i=0}^{d} E_{i}(A^{*}-a^{*}_{r} I)E_{r} \notag \\
                         &= E_{r}(A^{*}-a^{*}_{r} I)E_{r}+E_{s}A^{*}E_{r}. \notag
\end{align}
We have $E_{r}(A^{*}-a^{*}_{r} I)E_{r}=0$ by Proposition \ref{prop:EA*E}, so $(A^{*}-a^{*}_{r}I)E_{r}=E_{s}A^{*}E_{r}$. We apply both sides of this equation to $V$ and use Lemma \ref{lem:ijadjacent} to obtain $(A^{*}-a^{*}_{r} I)E_{r}V=E_{s}V$. Now take $\kappa=a^{*}_{r}$.

\medskip

\noindent (ii) $\Rightarrow$ (i). For $0\leq i\leq d$ such that $i\neq r$, we apply $E_{i}$ to both sides of the equation $(A^{*}-\kappa I)E_{r}V=E_{s}V$ and obtain $E_{i}A^{*}E_{r}V=\delta_{i,s}E_{i}V$. By Definition \ref{def:delta} and Lemma \ref{lem:ijadjacent}, vertex $r$ is adjacent to vertex $s$ and no other vertices.

\medskip

\noindent Suppose conditions (i) and (ii) hold. In the equation of (ii), apply $E_{r}$ to both sides to obtain $E_{r}(A^{*}-\kappa I)E_{r}V=0$. Therefore, $E_{r}(A^{*}-\kappa I)E_{r}=0$, so $E_{r}A^{*}E_{r}=\kappa E_{r}$. Now $\kappa=a^{*}_{r}$ in view of Proposition \ref{prop:EA*E}. \hfill $\Box$ \\

\noindent We record a result for later use.

\begin{lemma} \label{lem:A*=kI}
With reference to Assumption \ref{assum:A*}, the following {\rm (i)}--{\rm (iii)} are equivalent.
\begin{enumerate}
\item[\rm (i)] $A^{*}$ is a scalar multiple of $I$.
\item[\rm (ii)] The scalar $\theta^{*}_{i}$ is independent of $i$ for $0\leq i\leq d$.
\item[\rm (iii)] There exists an integer $r$ $(0\leq r\leq d)$ such that $\theta^{*}_{i}=a^{*}_{r}$ for $0\leq i\leq d$.
\end{enumerate}
Suppose conditions {\rm (i)}--{\rm (iii)} hold. Then $a^{*}_{i}$ is independent of $i$ for $0\leq i\leq d$.
\end{lemma}
\noindent {\it Proof:} Immediate from Assumption \ref{assum:A*} and Definition \ref{def:a*}. \hfill $\Box$ \\

\section{Recognizing leaves in $\Delta$ (part 2)} \label{sec:leaf_alg}

\begin{definition} \label{def:feasible}
\rm
With reference to Assumption \ref{assum:A*}, let $\{v_{i}\}_{i=0}^{d}$ denote a basis of $V$. We say that this basis is {\it feasible} whenever $v_{i}\in E^{*}_{i}V$ for $0\leq i\leq d$.
\end{definition}

\noindent With reference to Assumption \ref{assum:A*}, let $\{v_{i}\}_{i=0}^{d}$ denote a feasible basis of $V$. By Proposition \ref{prop:E*AE*} and Definition \ref{def:feasible}, the matrices representing $A$ and $A^{*}$ with respect to $\{v_{i}\}_{i=0}^{d}$ are
\begin{equation} \label{eq:A_A*_mat_gen}
A:\left(
\begin{array}
{ c c c c c c}
  a_{0} & b_{0} &       &       &       & {\bf 0} \\
  c_{1} & a_{1} & b_{1} &       &       & \\
        & c_{2} & \cdot & \cdot &       & \\
        &       & \cdot & \cdot & \cdot & \\
        &       &       & \cdot & \cdot & b_{d-1} \\
{\bf 0} &       &       &       & c_{d} & a_{d}
\end{array}
\right)
\qquad
A^{*}:\left(
\begin{array}
{ c c c c c c}
\theta^{*}_{0} &                &       &       &       & {\bf 0} \\
               & \theta^{*}_{1} &       &       &       & \\
               &                & \cdot &       &       & \\
               &                &       & \cdot &       & \\
               &                &       &       & \cdot & \\
{\bf 0}        &                &       &       &       & \theta^{*}_{d}
\end{array}
\right),
\end{equation}
where the scalars $\{a_{i}\}_{i=0}^{d}$ are from Definition \ref{def:a} and each of the scalars $\{b_{i}\}_{i=0}^{d-1}$ and $\{c_{i}\}_{i=1}^{d}$ is nonzero. For notational convenience, let $b_{d}=0$ and $c_{0}=0$.

\medskip

\noindent Observe that for $1\leq i\leq d$, the product $b_{i-1}c_{i}$ is independent of our choice of feasible basis. However, the scalars $\{b_{i}\}_{i=0}^{d-1}$ depend on the choice of feasible basis. Therefore, it is natural to ask what the possibilities are for $\{b_{i}\}_{i=0}^{d-1}$. The following lemma addresses this question.

\begin{lemma} \label{lem:A_A*_mat2}
With reference to Assumption \ref{assum:A*}, let $\{\beta_{i}\}_{i=0}^{d-1}$ denote a sequence of nonzero scalars taken from $\fld$. Then there exists a feasible basis of $V$ with respect to which the matrix representing $A$ has $(i,i+1)$-entry $\beta_{i}$ for $0\leq i\leq d-1$.
\end{lemma}
\noindent {\it Proof:} Let $\{v_{i}\}_{i=0}^{d}$ denote any feasible basis for $V$ and let $\{b_{i}\}_{i=0}^{d-1}$ denote the corresponding scalars from (\ref{eq:A_A*_mat_gen}). Define $v_{i}'=\frac{\beta_{0}\beta_{1}\cdots \beta_{i-1}}{b_{0}b_{1} \cdots b_{i-1}}v_{i}$ for $0\leq i\leq d$. Then $\{v_{i}'\}_{i=0}^{d}$ is a basis for $V$ which satisfies the requirements of the lemma. \hfill $\Box$ \\

\begin{definition}
\rm
With reference to Assumption \ref{assum:A*}, it follows by Lemma \ref{lem:A_A*_mat2} that there exists a feasible basis of $V$ such that $b_{i}=1$ for $0\leq i\leq d-1$. We call this basis the \emph{normalized feasible basis of $V$}.
\end{definition}

\noindent Let $\lambda$ denote an indeterminate. Let $\fld[\lambda]$ denote the $\fld$-algebra consisting of the polynomials in $\lambda$ that have all coefficients in $\fld$.

\begin{definition} \label{def:u_i}
\rm
With reference to Assumption \ref{assum:A*}, let $\{v_{i}\}_{i=0}^{d}$ denote a feasible basis for $V$. Define a sequence of polynomials $\{u_{i}\}_{i=0}^{d+1}$ in $\fld[\lambda]$ by
\begin{align}
u_{0} &= 1, \label{eq:u_0} \\
\lambda u_{i} &= c_{i}u_{i-1}+a_{i}u_{i}+b_{i}u_{i+1} \qquad \qquad (0 \leq i \leq d-1), \label{eq:u_i} \\
\lambda u_{d} &= c_{d}u_{d-1}+a_{d}u_{d}+\frac{u_{d+1}}{b_{0}b_{1}\cdots b_{d-1}}, \label{eq:u_d+1}
\end{align}
where $u_{-1}=0$. Observe that for $0\leq i\leq d+1$, the polynomial $u_{i}$ has degree $i$. Moreover, the coefficient of $\lambda^{i}$ in $u_{i}$ equals $(b_{0}b_{1}\cdots b_{i-1})^{-1}$ if $0\leq i\leq d$ and $1$ if $i=d+1$. We say that the sequence $\{u_{i}\}_{i=0}^{d+1}$ \emph{corresponds} to the feasible basis $\{v_{i}\}_{i=0}^{d}$. 
\end{definition}

\begin{definition} \label{def:p_i}
\rm
With reference to Assumption \ref{assum:A*}, let $\{p_{i}\}_{i=0}^{d+1}$ denote the polynomial sequence that corresponds to the normalized feasible basis of $V$. Observe that $p_{i}$ is monic of degree $i$ for $0\leq i\leq d+1$.
\end{definition}

\noindent We adopt the following assumption for the remainder of the section.

\begin{assumption} \label{assum:fix_feas}
\rm
With reference to Assumption \ref{assum:A*}, fix a feasible basis $\{v_{i}\}_{i=0}^{d}$ of $V$. Let $\{b_{i}\}_{i=0}^{d-1}$ and $\{c_{i}\}_{i=1}^{d}$ denote the corresponding scalars from (\ref{eq:A_A*_mat_gen}). Let $\{u_{i}\}_{i=0}^{d+1}$ denote the corresponding sequence of polynomials from Definition \ref{def:u_i}.
\end{assumption}

\noindent Recall the polynomials $\{p_{i}\}_{i=0}^{d+1}$ from Definition \ref{def:p_i}. From the perspective of Assumption \ref{assum:fix_feas}, these polynomials appear as follows.

\begin{lemma}
With reference to Assumptions \ref{assum:A*} and \ref{assum:fix_feas},
\begin{align}
p_{0} &= 1, \label{eq:p_0} \\
\lambda p_{i} &= b_{i-1}c_{i}p_{i-1}+a_{i}p_{i}+p_{i+1} \qquad \qquad (0 \leq i \leq d-1), \label{eq:p_i} \\
\lambda p_{d} &= b_{d-1}c_{d}p_{d-1}+a_{d}p_{d}+p_{d+1}, \label{eq:p_d+1}
\end{align}
where $p_{-1}=0$.
\end{lemma}
\noindent {\it Proof:} Consider the entries of the matrix representing $A$ with respect to the normalized feasible basis. For $1\leq i\leq d$, the product of the $(i-1,i)$-entry and the $(i,i-1)$-entry is $b_{i-1}c_{i}$. The $(i-1,i)$-entry is $1$, so the $(i,i-1)$-entry is $b_{i-1}c_{i}$. The result follows by Definitions \ref{def:u_i} and \ref{def:p_i}. \hfill $\Box$ \\

\begin{lemma} \label{lem:p_u}
With reference to Assumptions \ref{assum:A*} and \ref{assum:fix_feas},
\begin{align}
  u_{i} &= \frac{p_{i}}{b_{0}b_{1}\cdots b_{i-1}} \qquad \qquad (0\leq i\leq d), \notag \\
u_{d+1} &= p_{d+1}. \notag
\end{align}
\end{lemma}
\noindent {\it Proof:} Compare (\ref{eq:u_0})--(\ref{eq:u_d+1}) with (\ref{eq:p_0})--(\ref{eq:p_d+1}). \hfill $\Box$ \\

\begin{lemma} \label{lem:recur_u}
With reference to Assumptions \ref{assum:A*} and \ref{assum:fix_feas}, let $v$ denote a nonzero vector in $V$ and write $v=\sum_{i=0}^{d}\alpha_{i}v_{i}$. Let $\theta \in \fld$. Then the following {\rm (i)--(iii)} are equivalent.
\begin{enumerate}
\item[\rm (i)] The vector $v$ is an eigenvector for $A$ with eigenvalue $\theta$.
\item[\rm (ii)] For $0\leq i\leq d$,
\begin{equation} \label{eq:alpha_i_TTR}
c_{i}\alpha_{i-1}+a_{i}\alpha_{i}+b_{i}\alpha_{i+1}=\theta\alpha_{i},
\end{equation}
where $\alpha_{-1}$ and $\alpha_{d+1}$ are indeterminates.
\item[\rm (iii)] $\alpha_{i}=u_{i}(\theta)\alpha_{0}$ for $0\leq i\leq d$ and $u_{d+1}(\theta)=0$.
\end{enumerate}
Suppose conditions {\rm (i)--(iii)} hold. Then $\alpha_{0}\neq 0$.
\end{lemma}
\noindent {\it Proof:} By (\ref{eq:A_A*_mat_gen}),
\begin{equation} \label{eq:Av_i}
Av_{i}=b_{i-1}v_{i-1}+a_{i}v_{i}+c_{i+1}v_{i+1} \qquad (0\leq i\leq d),
\end{equation}
where $v_{-1}=0$ and $v_{d+1}=0$.

\medskip

\noindent (i) $\Rightarrow$ (ii). By assumption, $(A-\theta I)v=0$. In this equation, evaluate $v$ using $v=\sum_{i=0}^{d}\alpha_{i}v_{i}$, and simplify the result using (\ref{eq:Av_i}) and the fact that the $\{v_{i}\}_{i=0}^{d}$ are linearly independent. The result follows.

\medskip

\noindent (ii) $\Rightarrow$ (i). Using (\ref{eq:alpha_i_TTR}), along with $v=\sum_{i=0}^{d}\alpha_{i}v_{i}$ and (\ref{eq:Av_i}), we routinely obtain $Av=\theta v$.

\medskip

\noindent (ii) $\Rightarrow$ (iii). Let $\lambda =\theta$ in (\ref{eq:u_0})--(\ref{eq:u_d+1}) and compare the results to (\ref{eq:alpha_i_TTR}). This establishes $\alpha_{i}=u_{i}(\theta)\alpha_{0}$ for $0\leq i\leq d$ and $u_{d+1}(\theta)\alpha_{0}=0$. Note that $\alpha_{0}$ is nonzero; otherwise $\alpha_{i}=0$ for $0\leq i\leq d$, contradicting the assumption that $v$ is nonzero. Now $u_{d+1}(\theta)=0$.

\medskip

\noindent (iii) $\Rightarrow$ (ii). To verify (\ref{eq:alpha_i_TTR}), we eliminate $\alpha_{i-1}$, $\alpha_{i}$, $\alpha_{i+1}$ using $\alpha_{j}=u_{j}(\theta)\alpha_{0}$ ($0\leq j\leq d$) and evaluate the results using Definition \ref{def:u_i}.

\medskip

\noindent Suppose conditions (i)--(iii) hold. It was mentioned in the proof of (ii) $\Rightarrow$ (iii) that $\alpha_{0}\neq 0$. \hfill $\Box$ \\

\begin{corollary} \label{cor:u_d+1}
With reference to Assumptions \ref{assum:A*} and \ref{assum:fix_feas}, the polynomial $u_{d+1}$ is the characteristic polynomial for $A$.
\end{corollary}
\noindent {\it Proof:} By Definition \ref{def:u_i}, $u_{d+1}$ is monic with degree $d+1$. By Lemma \ref{lem:recur_u}, the $d+1$ eigenvalues of $A$ are all roots of $u_{d+1}$. The result follows. \hfill $\Box$ \\

\noindent With reference to Assumptions \ref{assum:A*} and \ref{assum:fix_feas}, let $\theta$ denote an eigenvalue of $A$. In Lemma \ref{lem:recur_u}(iii), we encountered the sequence $\{u_{i}(\theta)\}_{i=0}^{d}$. In the theory of distance-regular graphs, this sequence is called the cosine sequence for $\theta$. Motivated by this, we call the sequence $\{u_{i}(\theta)\}_{i=0}^{d}$ the \emph{cosine sequence for $\theta$ with respect to $\{v_{i}\}_{i=0}^{d}$}. Sometimes it is clear from the context what the basis $\{v_{i}\}_{i=0}^{d}$ is. In this case, we will refer to $\{u_{i}(\theta)\}_{i=0}^{d}$ as the \emph{cosine sequence for $\theta$}.

\begin{lemma} \label{lem:cos_recur}
With reference to Assumptions \ref{assum:A*} and \ref{assum:fix_feas}, let $\theta$ and $\{\alpha_{i}\}_{i=0}^{d}$ denote scalars in $\fld$. Then the following {\rm (i)}, {\rm (ii)} are equivalent.
\begin{enumerate}
\item[\rm (i)] The scalar $\theta$ is an eigenvalue for $A$ and $\{\alpha_{i}\}_{i=0}^{d}$ is the corresponding cosine sequence.
\item[\rm (ii)] $\alpha_{0}=1$ and
\begin{equation} \label{eq:cos_TTR}
c_{i}\alpha_{i-1}+a_{i}\alpha_{i}+b_{i}\alpha_{i+1}=\theta\alpha_{i} \qquad (0 \leq i \leq d),
\end{equation}
where $\alpha_{-1}$ and $\alpha_{d+1}$ are indeterminates.
\end{enumerate}
Suppose conditions {\rm (i)} and {\rm (ii)} hold. Then $\sum_{i=0}^{d}\alpha_{i}v_{i}$ is an eigenvector for $A$ with eigenvalue $\theta$.
\end{lemma}
\noindent {\it Proof:} This is an immediate consequence of Lemma \ref{lem:recur_u}. \hfill $\Box$ \\

\noindent To motivate the upcoming results, we have some comments. Choose an integer $r$ $(0\leq r\leq d)$. Let $\{\alpha_{i}\}_{i=0}^{d}$ denote the cosine sequence for $\theta_{r}$, so that
\begin{equation} \label{eq:alpha}
\alpha_{i}=u_{i}(\theta_{r}) \qquad \qquad (0\leq i\leq d).
\end{equation}
Define $v=\sum_{i=0}^{d}\alpha_{i}v_{i}$. The vector $v$ is nonzero by construction and contained in $E_{r}V$ by Lemma \ref{lem:cos_recur}. Let $w=(A^{*}-a^{*}_{r}I)v$. Recall $A^{*}v_{i}=\theta^{*}_{i}v_{i}$ for $0\leq i\leq d$. Therefore, $w=\sum_{i=0}^{d}\alpha'_{i}v_{i}$, where
\begin{equation} \label{eq:alpha'}
\alpha'_{i}=(\theta^{*}_{i}-a^{*}_{r})\alpha_{i} \qquad \qquad (0\leq i\leq d).
\end{equation}

\begin{lemma} \label{lem:ar*neqtheta*0}
With reference to Assumptions \ref{assum:A*} and \ref{assum:fix_feas}, fix an integer $r$ such that $0\leq r\leq d$. In the diagram $\Delta$, assume vertex $r$ is adjacent to exactly one other vertex. Then $a^{*}_{r}\neq \theta^{*}_{0}$.
\end{lemma}
\noindent {\it Proof:} We refer to the vector $w$ from the paragraph preceding this lemma. Let $s$ denote the vertex in $\Delta$ that is adjacent to $r$. We have $(A^{*}-a^{*}_{r}I)E_{r}V=E_{s}V$ by Lemma \ref{lem:A*-a*I}, so $0\neq w\in E_{s}V$. In the sum $w=\sum_{i=0}^{d}\alpha'_{i}v_{i}$, the coefficient $\alpha'_{0}$ is nonzero by the final assertion of Lemma \ref{lem:recur_u}. By this and (\ref{eq:alpha'}), we find $a^{*}_{r}\neq \theta^{*}_{0}$. \hfill $\Box$ \\

\begin{proposition} \label{prop:genstailTTR}
With reference to Assumptions \ref{assum:A*} and \ref{assum:fix_feas}, fix distinct integers $r$ and $s$ such that $0\leq r, s\leq d$. Then the following {\rm (i)}, {\rm (ii)} are equivalent.
\begin{enumerate}
\item[\rm (i)] In the diagram $\Delta$, vertex $r$ is adjacent to vertex $s$ and no other vertices.
\item[\rm (ii)] The cosine sequence $\{\alpha_{i}\}_{i=0}^{d}$ for $\theta_{r}$ satisfies
\begin{equation} \label{eq:genstailTTR}
c_{i}\theta^{*}_{i-1}\alpha_{i-1}+a_{i}\theta^{*}_{i}\alpha_{i}+b_{i}\theta^{*}_{i+1}\alpha_{i+1}-\theta_{r}\theta^{*}_{i}\alpha_{i}=(\theta_{s}-\theta_{r})(\theta^{*}_{i}-a^{*}_{r})\alpha_{i} \quad (0\leq i\leq d),
\end{equation}
where each of $\alpha_{-1}$, $\alpha_{d+1}$, $\theta^{*}_{-1}$, and $\theta^{*}_{d+1}$ is indeterminate. Furthermore, there exists an integer $i$ ($0\leq i\leq d$) such that the right-hand side of {\rm (\ref{eq:genstailTTR})} is not equal to $0$.
\end{enumerate}
\end{proposition}
\noindent {\it Proof:} We refer to the vectors $v,w$ from the paragraph preceding Lemma \ref{lem:ar*neqtheta*0}.

\medskip

\noindent (i) $\Rightarrow$ (ii). The vector $w$ spans $E_{s}V$ by Lemma \ref{lem:A*-a*I}, so $Aw=\theta_{s}w$ and $w$ is nonzero. Applying Lemma \ref{lem:recur_u} to $w$, 
\begin{equation}
c_{i}\alpha'_{i-1}+a_{i}\alpha'_{i}+b_{i}\alpha'_{i+1}=\theta_{s}\alpha'_{i} \qquad \qquad (0\leq i\leq d), \notag
\end{equation}
where $\alpha'_{-1}$ and $\alpha'_{d+1}$ are indeterminates. Evaluating this equation first using (\ref{eq:alpha'}) and then using (\ref{eq:alpha_i_TTR}), we obtain (\ref{eq:genstailTTR}). By construction, $\theta_{s}\neq \theta_{r}$. Also, $w$ is nonzero, so $\{\alpha'_{i}\}_{i=0}^{d}$ are not all zero. By these comments and (\ref{eq:alpha'}), there exists an integer $i$ ($0\leq i\leq d$) such that the right-hand side of (\ref{eq:genstailTTR}) is not equal to $0$.

\medskip

\noindent (ii) $\Rightarrow$ (i). We show that $(A^{*}-a^{*}_{r} I)E_{r}V=E_{s}V$. We mentioned earlier that $0\neq v\in E_{r}V$. Using (\ref{eq:A_A*_mat_gen}), rewrite (\ref{eq:genstailTTR}) as $AA^{*}v=\theta_{s}(A^{*}-a^{*}_{r} I)v+a^{*}_{r} \theta_{r}v$. In this equation, we rearrange terms and use $Av=\theta_{r}v$ to obtain $Aw=\theta_{s}w$. Therefore, $w\in E_{s}V$. By assumption, there exists an integer $i$ ($0\leq i\leq d$) such that $(\theta^{*}_{i}-a^{*}_{r})\alpha_{i}\neq 0$. By this and (\ref{eq:alpha'}), $w\neq 0$. By these comments, $(A^{*}-a^{*}_{r} I)E_{r}V=E_{s}V$. The result follows in view of Lemma \ref{lem:A*-a*I}. \hfill $\Box$ \\

\begin{proposition} \label{prop:AWstail_gen}
With reference to Assumptions \ref{assum:A*} and \ref{assum:fix_feas}, fix distinct integers $r$ and $s$ such that $0\leq r,s\leq d$. Then the following {\rm (i)}, {\rm (ii)} are equivalent.
\begin{enumerate}
\item[\rm (i)] In the diagram $\Delta$, vertex $r$ is adjacent to vertex $s$ and no other vertices.
\item[\rm (ii)] $a^{*}_{r} \neq \theta^{*}_{0}$ and
\begin{equation} \label{eq:AWD_gen}
u_{i}(\theta_{s})=u_{i}(\theta_{r})\frac{\theta^{*}_{i}-a^{*}_{r}}{\theta^{*}_{0}-a^{*}_{r}} \qquad \qquad (0\leq i\leq d),
\end{equation}
where the polynomials $\{u_{i}\}_{i=0}^{d}$ are from Assumption \ref{assum:fix_feas}.
\end{enumerate}
\end{proposition}
\noindent {\it Proof:} We refer to the vectors $v,w$ from the paragraph preceding Lemma \ref{lem:ar*neqtheta*0}.

\medskip

\noindent (i) $\Rightarrow$ (ii). Observe that $a^{*}_{r}\neq \theta^{*}_{0}$ by Lemma \ref{lem:ar*neqtheta*0}. By Lemma \ref{lem:A*-a*I}, $w$ is a nonzero vector in $E_{s}V$. By construction, $w$ is an eigenvector for $A$ with eigenvalue $\theta_{s}$. Applying Lemma \ref{lem:recur_u} to $w$, 
\begin{equation} \label{eq:alpha'_gen}
\alpha'_{i}=u_{i}(\theta_{s})\alpha'_{0} \qquad \qquad (0\leq i\leq d).
\end{equation}
To obtain (\ref{eq:AWD_gen}), evaluate each side of (\ref{eq:alpha'_gen}) using (\ref{eq:alpha'}) and simplify the result using (\ref{eq:alpha}).

\medskip

\noindent (ii) $\Rightarrow$ (i). We show $(A^{*}-a^{*}_{r}I)E_{r}V=E_{s}V$. Note that $\alpha'_{0}=\theta^{*}_{0}-a^{*}_{r}$ and this is nonzero by assumption. Evaluating (\ref{eq:alpha'}) using this fact and (\ref{eq:AWD_gen}), we obtain $\alpha'_{i}=u_{i}(\theta_{s})\alpha'_{0}$ for $0\leq i\leq d$. By Corollary \ref{cor:u_d+1}, $u_{d+1}(\theta_{s})=0$. The vector $w$ is nonzero by construction and contained in $E_{s}V$ by Lemma \ref{lem:recur_u}. We mentioned earlier that $0\neq v\in E_{r}V$. Therefore, $(A^{*}-a^{*}_{r}I)E_{r}V=E_{s}V$. The result follows in view of Lemma \ref{lem:A*-a*I}. \hfill $\Box$ \\

\noindent In many applications where Assumption \ref{assum:A*} is relevant, the matrix representing $A$ in (\ref{eq:A_A*_mat_gen}) has constant row sum. In the next section, we will adopt this assumption and investigate its consequences. In that investigation, the following results will be helpful.

\begin{proposition} \label{prop:rowsum}
With reference to Assumptions \ref{assum:A*} and \ref{assum:fix_feas}, the following {\rm (i)}, {\rm (ii)} are equivalent for $\theta \in \fld$.
\begin{enumerate}
\item[\rm (i)] The matrix on the left in {\rm (\ref{eq:A_A*_mat_gen})} has constant row sum $\theta$.
\item[\rm (ii)] The scalar $\theta$ is an eigenvalue of $A$ and $u_{i}(\theta)=1$ for $0\leq i\leq d$, where the polynomials $\{u_{i}\}_{i=0}^{d}$ are from Assumption \ref{assum:fix_feas}.
\end{enumerate}
\end{proposition}
\noindent {\it Proof:} (i) $\Rightarrow$ (ii). By assumption,
\begin{equation} \label{eq:rowsum=theta}
c_{i}+a_{i}+b_{i}=\theta \qquad \qquad (0\leq i\leq d),
\end{equation}
where $c_{0}=0$ and $b_{d}=0$. Let $v=\sum_{i=0}^{d}v_{i}$. By (\ref{eq:A_A*_mat_gen}) and (\ref{eq:rowsum=theta}), $Av=\theta v$, so $\theta$ is an eigenvalue of $A$. Let $\{\alpha_{i}\}_{i=0}^{d}$ denote the cosine sequence corresponding to $\theta$.  By Lemma \ref{lem:cos_recur}, $\alpha_{0}=1$ and (\ref{eq:cos_TTR}) holds. Comparing (\ref{eq:rowsum=theta}) to (\ref{eq:cos_TTR}), we find that $\alpha_{i}=1$ for $0\leq i\leq d$. In other words, $u_{i}(\theta)=1$ for $0\leq i\leq d$.

\medskip

\noindent (ii) $\Rightarrow$ (i). The cosine sequence corresponding to $\theta$ is $(1, 1, \ldots, 1)$. The result follows by Lemma \ref{lem:cos_recur}. \hfill $\Box$ \\

\begin{proposition} \label{prop:rowsum2}
With reference to Assumptions \ref{assum:A*} and \ref{assum:fix_feas}, the following {\rm (i)}, {\rm (ii)} are equivalent for $\theta \in \fld$.
\begin{enumerate}
\item[\rm (i)] There exists a feasible basis for $V$ with respect to which the matrix on the left in {\rm (\ref{eq:A_A*_mat_gen})} has constant row sum $\theta$.
\item[\rm (ii)] The scalar $\theta$ is an eigenvalue of $A$ and $u_{i}(\theta)\neq 0$ for $0\leq i\leq d$, where the polynomials $\{u_{i}\}_{i=0}^{d}$ are from Assumption \ref{assum:fix_feas}.
\end{enumerate}
\end{proposition}
\noindent {\it Proof:} (i) $\Rightarrow$ (ii). Let $\{v'_{i}\}_{i=0}^{d}$ denote a feasible basis for $V$ with respect to which the matrix on the left in (\ref{eq:A_A*_mat_gen}) has constant row sum $\theta$. Let $\{b'_{i}\}_{i=0}^{d-1}$ denote the corresponding scalars from (\ref{eq:A_A*_mat_gen}), and let $\{u'_{i}\}_{i=0}^{d+1}$ denote the corresponding polynomial sequence from Definition \ref{def:u_i}. By Proposition \ref{prop:rowsum}, $\theta$ is an eigenvalue of $A$ and $u'_{i}(\theta)=1$ for $0\leq i\leq d$. By Lemma \ref{lem:p_u}, $p_{i}(\theta)=u_{i}(\theta)b_{0}b_{1}\cdots b_{i-1}$ and $p_{i}(\theta)=u'_{i}(\theta)b'_{0}b'_{1}\cdots b'_{i-1}$ for $0\leq i\leq d$. We have $u'_{i}(\theta)\neq 0$ for $0\leq i\leq d$. Also, $\{b'_{i}\}_{i=0}^{d-1}$ are all nonzero. Therefore, $p_{i}(\theta)\neq 0$ for $0\leq i\leq d$. Consequently, $u_{i}(\theta)\neq 0$ for $0\leq i\leq d$.

\medskip

\noindent (ii) $\Rightarrow$ (i). For notational convenience, define $\gamma_{i}=\frac{u_{i}(\theta)}{u_{i-1}(\theta)}$ for $1\leq i\leq d$. By assumption, the scalars $\{\gamma_{i}\}_{i=1}^{d}$ are nonzero. By Lemma \ref{lem:A_A*_mat2}, there exists a feasible basis of $V$ with respect to which the matrix representing $A$ has $(i-1,i)$-entry $b_{i-1}\gamma_{i}$ for $1\leq i\leq d$. We have some comments regarding the entries of this matrix. For $1\leq i\leq d$, the product of the $(i-1,i)$-entry and the $(i,i-1)$-entry is $b_{i-1}c_{i}$. The $(i-1,i)$-entry is $b_{i-1}\gamma_{i}$, so the $(i,i-1)$-entry is $c_{i}\gamma^{-1}_{i}$. In row $0$, the entries that are potentially nonzero are the $(0,0)$-entry and the $(0,1)$-entry. These entries are $a_{0}$ and $b_{0}\gamma_{1}$, respectively. Therefore, the sum of the entries in row $0$ is $a_{0}+b_{0}\gamma_{1}$. For $1\leq i\leq d-1$, the entries in row $i$ that are potentially nonzero are the $(i,i-1)$-entry, the $(i,i)$-entry, and the $(i,i+1)$-entry. These entries are $c_{i}\gamma^{-1}_{i}$, $a_{i}$, and $b_{i}\gamma_{i+1}$, respectively. Therefore, the sum of the entries in the $i^{th}$ row is $c_{i}\gamma^{-1}_{i}+a_{i}+b_{i}\gamma_{i+1}$. In row $d$, the entries that are potentially nonzero are the $(d,d-1)$-entry and the $(d,d)$-entry. These entries are $c_{d}\gamma^{-1}_{d}$ and $a_{d}$, respectively. Therefore, the sum of the entries in row $d$ is $c_{d}\gamma^{-1}_{d}+a_{d}$. For each row, the sum of the entries is equal to $\theta$ by Definition \ref{def:u_i} and Corollary \ref{cor:u_d+1}. \hfill $\Box$ \\

\section{Recognizing leaves in $\Delta$ (part 3)} \label{sec:leaf3}

\begin{assumption} \label{assum:CRS}
\rm
With reference to Assumption \ref{assum:A*}, let $\theta_{r}$ denote an eigenvalue for $A$. Fix a feasible basis $\{v_{i}\}_{i=0}^{d}$ of $V$. Assume the matrix representing $A$ with respect to this basis has constant row sum $\theta_{r}$. We retain the notation from (\ref{eq:A_A*_mat_gen}) for the matrices representing $A$ and $A^{*}$.
\end{assumption}

\noindent With reference to Assumptions \ref{assum:A*} and \ref{assum:CRS}, note that each entry in the cosine sequence for $\theta_{r}$ is equal to $1$ by Proposition \ref{prop:rowsum}. Let $v=\sum_{i=0}^{d}v_{i}$ and note that $v$ is nonzero. Also, $Av=\theta_{r}v$ by (\ref{eq:A_A*_mat_gen}) and (\ref{eq:rowsum=theta}). Therefore, $v$ spans $E_{r}V$.

\begin{proposition}
With reference to Assumptions \ref{assum:A*} and \ref{assum:CRS}, the following {\rm (i)}, {\rm (ii)} are equivalent for $0\leq s\leq d$ with $r\neq s$.
\begin{enumerate}
\item[\rm (i)] In the diagram $\Delta$, vertex $r$ is adjacent to vertex $s$ and no other vertices.
\item[\rm (ii)] For $0\leq i\leq d$,
\begin{equation} \label{eq:stailTTR}
c_{i}\theta^{*}_{i-1}+a_{i}\theta^{*}_{i}+b_{i}\theta^{*}_{i+1}-\theta_{r}\theta^{*}_{i}=(\theta_{s}-\theta_{r})(\theta^{*}_{i}-a^{*}_{r}),
\end{equation}
where $\theta^{*}_{-1}$ and $\theta^{*}_{d+1}$ are indeterminates. Furthermore, it is not the case that $A^{*}$ is a scalar multiple of the identity.
\end{enumerate}
\end{proposition}
\noindent {\it Proof:} By Proposition \ref{prop:genstailTTR} and the discussion below Assumption \ref{assum:CRS}, condition (i) holds if and only if both (\ref{eq:stailTTR}) holds for $0\leq i\leq d$ and there exists an integer $i$ ($0\leq i\leq d$) such that the right-hand side of (\ref{eq:stailTTR}) is not equal to $0$. However, the right-hand side of (\ref{eq:stailTTR}) equals $0$ for $0\leq i\leq d$ if and only if $\theta^{*}_{i}=a^{*}_{r}$ for $0\leq i\leq d$. The result follows by Lemma \ref{lem:A*=kI}. \hfill $\Box$ \\

\begin{proposition}
With reference to Assumptions \ref{assum:A*} and \ref{assum:CRS}, for $0\leq s\leq d$ such that $r\neq s$, the following {\rm (i)}, {\rm (ii)} are equivalent.
\begin{enumerate}
\item[\rm (i)] In the diagram $\Delta$, vertex $r$ is adjacent to vertex $s$ and no other vertices.
\item[\rm (ii)] $a^{*}_{r}\neq \theta^{*}_{0}$ and
\begin{equation}
u_{i}(\theta_{s})=\frac{\theta^{*}_{i}-a^{*}_{r}}{\theta^{*}_{0}-a^{*}_{r}} \qquad \qquad (0\leq i\leq d), \notag
\end{equation}
where the polynomials $\{u_{i}\}_{i=0}^{d}$ are from Assumption \ref{assum:fix_feas}.
\end{enumerate}
\end{proposition}
\noindent {\it Proof:} The result follows by Proposition \ref{prop:AWstail_gen} and the discussion below Assumption \ref{assum:CRS}. \hfill $\Box$ \\

\section{Appendix A: An algorithm for recognizing a leaf}

Given the conditions of Assumption \ref{assum:A*}, define a diagram $\Delta$ as in Definition \ref{def:delta}. We present an algorithm designed to recognize a leaf in $\Delta$. This algorithm is based on the theory developed in Sections \ref{sec:a*} and \ref{sec:leaf_alg}. Suppose we are given two distinct integers $r$ and $s$ such that $0\leq r,s\leq d$. Our algorithm checks whether vertex $r$ is adjacent to vertex $s$ and no other vertices. We say that the ordered pair $(r, s)$ is \emph{confirmed} whenever this is the case. Otherwise, $(r, s)$ is said to be \emph{denied}. To avoid trivialities, we assume $d\geq 2$, so that $\Delta$ contains at least three vertices. We also assume that $\theta^{*}_{i}\neq \theta^{*}_{0}$ for $1\leq i \leq d$, so that $\Delta$ is connected by Proposition \ref{prop:deltaconnected}. Fix a feasible basis $\{v_{i}\}_{i=0}^{d}$ for $V$. As we saw in (\ref{eq:A_A*_mat_gen}), the matrices representing $A$ and $A^{*}$ with respect to this basis are:
\begin{equation}
A:\left(
\begin{array}
{ c c c c c c}
  a_{0} & b_{0} &       &       &       & {\bf 0} \\
  c_{1} & a_{1} & b_{1} &       &       & \\
        & c_{2} & \cdot & \cdot &       & \\
        &       & \cdot & \cdot & \cdot & \\
        &       &       & \cdot & \cdot & b_{d-1} \\
{\bf 0} &       &       &       & c_{d} & a_{d}
\end{array}
\right)
\qquad
A^{*}:\left(
\begin{array}
{ c c c c c c}
\theta^{*}_{0} &                &       &       &       & {\bf 0} \\
               & \theta^{*}_{1} &       &       &       & \\
               &                & \cdot &       &       & \\
               &                &       & \cdot &       & \\
               &                &       &       & \cdot & \\
{\bf 0}        &                &       &       &       & \theta^{*}_{d}
\end{array}
\right), \notag
\end{equation}
where the scalars $\{a_{i}\}_{i=0}^{d}$ are from Definition \ref{def:a} and $b_{i-1}c_{i}\neq 0$ for $1\leq i\leq d$. Let 
\begin{equation}
\alpha_{0}=1, \qquad \alpha_{1}=\frac{\theta_{r}-a_{0}}{b_{0}}, \qquad \kappa=\frac{\theta^{*}_{1}(\theta_{r}-a_{0})-\theta^{*}_{0}(\theta_{s}-a_{0})}{\theta_{r}-\theta_{s}}. \notag
\end{equation}
The algorithm consists of the following steps (i)--(v):

\begin{enumerate}
\item Let $j=1$.
\item Calculate $\alpha_{j+1}$ using $c_{j}\alpha_{j-1}+a_{j}\alpha_{j}+b_{j}\alpha_{j+1}=\theta_{r}\alpha_{j}$.
\item Check if $c_{j}\theta^{*}_{j-1}\alpha_{j-1}+a_{j}\theta^{*}_{j}\alpha_{j}+b_{j}\theta^{*}_{j+1}\alpha_{j+1}=(\theta_{s}\theta^{*}_{j}+\kappa(\theta_{r}-\theta_{s}))\alpha_{j}$. If so, then proceed to the next step. If not, then $(r, s)$ is denied.
\item Add $1$ to $j$.
\item Check if $j\leq d-1$. If so, then go back to step (ii). If not, then $(r, s)$ is confirmed.
\end{enumerate}

\section{Appendix B: An algorithm for recognizing a leaf assuming $A$ has constant row sum}

Given the conditions of Assumption \ref{assum:A*}, define a diagram $\Delta$ as in Definition \ref{def:delta}. In this appendix, we present a second algorithm designed to recognize a leaf in $\Delta$. Our setup is the same as in Appendix A with the additional assumption that, with respect to a fixed feasible basis for $V$, the matrix representing $A$ has constant row sum $\theta_{r}$. With respect to this basis, the matrices representing $A$ and $A^{*}$ are:
\begin{equation}
A:\left(
\begin{array}
{ c c c c c c}
  a_{0} & b_{0} &       &       &       & {\bf 0} \\
  c_{1} & a_{1} & b_{1} &       &       & \\
        & c_{2} & \cdot & \cdot &       & \\
        &       & \cdot & \cdot & \cdot & \\
        &       &       & \cdot & \cdot & b_{d-1} \\
{\bf 0} &       &       &       & c_{d} & a_{d}
\end{array}
\right)
\qquad
A^{*}:\left(
\begin{array}
{ c c c c c c}
\theta^{*}_{0} &                &       &       &       & {\bf 0} \\
               & \theta^{*}_{1} &       &       &       & \\
               &                & \cdot &       &       & \\
               &                &       & \cdot &       & \\
               &                &       &       & \cdot & \\
{\bf 0}        &                &       &       &       & \theta^{*}_{d}
\end{array}
\right), \notag
\end{equation}
where the scalars $\{a_{i}\}_{i=0}^{d}$ are from Definition \ref{def:a}, $b_{i-1}c_{i}\neq 0$ for $1\leq i\leq d$, and $c_{i}+a_{i}+b_{i}=\theta_{r}$ for $0\leq i\leq d$. Let 
\begin{equation}
\kappa=\frac{\theta^{*}_{1}b_{0}-\theta^{*}_{0}(\theta_{s}-a_{0})}{\theta_{r}-\theta_{s}}. \notag
\end{equation}
The algorithm consists of the following steps (i)--(iv):

\begin{enumerate}
\item Let $j=1$.
\item Check if $c_{j}\theta^{*}_{j-1}+a_{j}\theta^{*}_{j}+b_{j}\theta^{*}_{j+1}=\theta_{s}\theta^{*}_{j}+\kappa(\theta_{r}-\theta_{s})$. If so, then proceed to the next step. If not, then $(r, s)$ is denied.
\item Add $1$ to $j$.
\item Check if $j\leq d-1$. If so, then go back to step (ii). If not, then $(r, s)$ is confirmed.
\end{enumerate}

\section{Acknowledgment}

This paper was written while the author was a graduate student at the University of Wisconsin-Madison. The author thanks his advisor, Paul Terwilliger, for offering many valuable ideas and suggestions.

\bigskip

\noindent Edward Hanson \hfil\break
\noindent Department of Mathematics \hfil\break
\noindent University of Wisconsin \hfil\break
\noindent 480 Lincoln Drive \hfil\break
\noindent Madison, WI 53706-1388 USA \hfil\break
\noindent email: {\tt hanson@math.wisc.edu }\hfil\break

\end{document}